\documentclass{amsart}
\usepackage{verbatim}

\newtheorem{thm}{Theorem}
\newtheorem*{tpl}{Triple Point Lemma}
\newtheorem*{incidence}{Incidence Formula}
\newtheorem{conjecture}{Conjecture}
\theoremstyle{definition} \newtheorem{define}{Definition}
\theoremstyle{remark}  \newtheorem*{example}{Example}

\newcommand{\prob}{\operatorname{P}}
\newcommand{\Var}{\operatorname{Var}}
\newcommand{\EE}{\mathbb{E}}

\newcommand{\FF}{\mathbb{F}_q}
\newcommand{\RR}{\mathbb{R}}

\newcommand{\primesum}{\sum\nolimits'}

\title[Kakeya problem]{On the finite field Kakeya problem \\ in two dimensions}

\author{X.W.C. Faber}
\address{X.W.C. Faber\\
Department of Mathematics\\
Columbia University\\
New York, NY  10027\\
USA} \email{xander@math.columbia.edu}



\date{October 18, 2005}

\begin{document}

\begin{abstract}
  A two-dimensional Besicovitch set over a finite field is a subset of
  the finite plane containing a line in each direction. In this paper,
  we conjecture a sharp lower bound for the size of such a subset and
  prove some results toward this conjecture.
\end{abstract}

\maketitle

\section{Introduction}

The classical Kakeya problem, posed in 1917 by Kakeya, asks for a compact region in
$\RR^2$ of minimum Lebesgue measure in which one can continuously turn a unit length
segment through a full $360^{\circ}$ rotation. By 1928 Besicovitch had proved that such a
region exists with arbitrarily small Lebesgue measure. Prior to this result, he also
constructed a compact subset of \textit{zero Lebesgue measure} containing a unit length
segment in any direction. (Of course, one can't continuously turn the segment in this
set.)

The finite field Kakeya problem, originally posited by Wolff in
\cite{Wo}, asks for the smallest subset of $\FF^n$ that contains a
line in each direction, where $\FF$ denotes the finite field with $q$
elements. A subset containing a line in each direction is called a
\textit{Besicovitch set}.  Wolff conjectured that there is a positive
constant $C = C(n)$ such that $\#B \geq Cq^n$ for any Besicovitch set
$B \subset \FF^n$. For $n=2$ he immediately proved that $\#B \geq
q^2/2$; Wolff's method actually gives $\#B \geq q(q+1)/2$. The finite
field Kakeya problem has also been investigated in \cite{Ro},
\cite{MT}, \cite{BKT}, and \cite{Ta}. These authors have concentrated
their efforts toward obtaining satisfactory asymptotic lower bounds
for a Besicovitch set in $\FF^n$ for $n \geq 3$.

In this paper, we focus our attention exclusively on Besicovitch sets
in $\FF^2$ and sharpen Wolff's lower bound by combinatorial methods.
The next section will be devoted to explaining new results. All of the
proofs will be given in section~\ref{Proof Section}.

It should also be noted that the recent work \cite{Co} builds upon the
techniques in the present article in order to improve one of the
results. See the remark following the statement of Theorem~\ref{Unconditional Theorem}.

\section{Results}

In $\FF^2$ a line is the set of solutions of an equation $ax + by = c$ with $a,b,c \in
\FF$. Write $\ell(m,b)$ and $\ell(\infty,a)$ for the lines $y=mx+b$ and $x=a$,
respectively.

\begin{define}
A \textit{Besicovitch set in $\FF^2$} is a set of points $B \subset \FF^2$ such that for
each $i \in \FF \cup \{\infty\}$ there exists $b_i \in \FF$ so that $\ell(i, b_i) \subset
B$.
\end{define}

The smallest Besicovitch sets will be those that are a union of lines with distinct
slopes. Regarding the size of such a set, we have the following:

\begin{incidence} \label{Incidence Proposition}
Suppose $B$ is a Besicovitch set with $B=\bigcup_{i \in \FF\cup\{\infty\}} \ell(i,b_i)$.
For $P \in \FF^2$, let $m_P$ be the number of these lines passing through $P$. Then
\begin{equation*}
\#B = \frac{q(q+1)}{2} + \sum_{P \in B} \frac{(m_P-1)(m_P-2)}{2}.
\end{equation*}
In particular, every Besicovitch set contains at least $\frac{q(q+1)}{2}$
points.
\end{incidence}
Our approach will be to study the intersections of lines in a Besicovitch set in order to
show that the sum in the Incidence Formula cannot be too small.

\begin{example} Consider the set
\begin{equation*} \label{Minimal Set}
B_0 = \left( \bigcup_{i \in \FF} \ell(i,-i^2) \right) \cup \ell(\infty,0).
\end{equation*}
One can calculate that
\begin{equation*} \label{Minimal Set Size}
\sum_{P \in B_0} \frac{(m_P-1)(m_P-2)}{2} =\begin{cases} 0 & \text{if $q$ is even},
\\\frac{q-1}{2} & \text{if $q$ is odd}. \end{cases}
\end{equation*}
We will perform this calculation in section~\ref{Proof Section}. If
$q$ is even, the set $B_0$ achieves the minimum cardinality allowed by
the Incidence Formula.  When $q$ is odd one might guess
that $\#B_0$ gives a sharp lower bound as well.
\end{example}

\begin{conjecture} \label{Size Conjecture}
If $q$ is odd, a Besicovitch set $B \subset \FF$ must have
\begin{equation*}
\sum_{P \in B} \frac{(m_P-1)(m_P-2)}{2} \geq \frac{q-1}{2}.
\end{equation*}
That is,
\begin{equation*}
\#B \geq \frac{q(q+1)}{2}+\frac{q-1}{2}.
\end{equation*}
\end{conjecture}

Our first main result is an improvement on the trivial lower bound
given by the Incidence Formula.

\begin{thm} \label{Unconditional Theorem}
Assume $q$ is odd. For any Besicovitch set $B$, we have
\begin{equation*} \label{Unconditional Estimate}
\sum_{P \in B} \frac{(m_P-1)(m_P-2)}{2} \geq \frac{q}{3}.
\end{equation*}
\end{thm}
Observe that this estimate immediately implies Conjecture~\ref{Size
  Conjecture} for $q=3, 5, 7$.

In \cite{Co}, Cooper has refined the strategy for proving Theorem~\ref{Unconditional
Theorem} and is able to obtain the stronger lower bound $(5q-1)/14$.

On the other hand, Theorem~\ref{Conditional Theorem} will give a sharp
conditional form of Conjecture~\ref{Size Conjecture}. Let us take a
moment to motivate the hypothesis of the theorem before we state it.

\begin{define} Let $q$ be odd. A Besicovitch set $B \subset \FF$ will
  be called \textit{small} if $\#B$ does not exceed the lower bound in
  Conjecture~\ref{Size Conjecture}.
\end{define}

The proof of the Incidence Formula will yield the same result if $B$ is a union of only
$q$ lines (as opposed to $q+1$). Therefore any subset of $\FF^2$ that is a union of $q$
lines in distinct directions must contain at least $q(q+1)/2$ points, with equality if
and only if no three of the lines share a common point. If we consider the set $B_0' =
\bigcup_{i \in \FF} \ell(i,-i^2)$, neglecting the vertical line in the above example, one
can see that $\#B_0' = q(q+1)/2$. Thus $B_0'$ has minimum cardinality among all sets
consisting of the union of $q$ lines in distinct directions. It seems plausible that such
a set has the best chance of yielding a small Besicovitch set when we adjoin one more
line. Note that any set constructed in this way will have all of its points of
multiplicity three lying on one line---the final line adjoined to the set.  Indeed, we
can prove that a Besicovitch set with this last property satisfies Conjecture~\ref{Size
  Conjecture}:

\begin{thm} \label{Conditional Theorem}
Assume $q$ is odd. Let $B=\bigcup_{i \in \FF\cup\{\infty\}} \ell(i,b_i)$ be a Besicovitch
set. Suppose there is $j \in \FF \cup \{\infty\}$ such that every point $P \in B$ with
$m_P \geq 3$ lies on the line $\ell(j,b_j)$. Then
\begin{equation*}
\sum_{P \in B} \frac{(m_P-1)(m_P-2)}{2} \geq \frac{q-1}{2}.
\end{equation*}
Equality holds if and only there are $(q-1)/2$ points $P \in B$ with $m_P=3$ and
no points with $m_P > 3$.
\end{thm}

It seems natural to state another conjecture in light of the above discussion.

\begin{conjecture} \label{Multiplicity Conjecture}
If $q$ is odd and $B=\bigcup_{i \in \FF \cup \{\infty\}} \ell(i,b_i)$ is a small Besicovitch
set, then there is $j \in \FF \cup \{\infty\}$ such that every point $P \in B$ with $m_P
\geq 3$ lies on the line $\ell(j,b_j)$.
\end{conjecture}

Now that we have two conjectures it seems reasonable to think about testing them via
computer calculation. By checking every Besicovitch set that is a union of $q+1$ lines,
we have learned that Conjectures~\ref{Size Conjecture}~and~\ref{Multiplicity Conjecture}
hold for $q\leq 13$ odd. Unfortunately, I haven't been able to construct an
algorithm for finding small Besicovitch sets that requires any fewer than about $O(q^q)$
steps. In order to try to disprove Conjecture~\ref{Size Conjecture}, one might randomly
select a collection of lines with distinct slopes and hope that it will yield a small
Besicovitch set. The following theorem shows that one is unlikely to get so lucky.

\begin{thm} \label{Probability Theorem}

\begin{enumerate}
\item[] \hspace{1in}

\item[(a)] The expected cardinality of a Besicovitch set formed by the union of $q+1$ randomly chosen lines
with distinct slopes is
\begin{equation*}
\left(1-\left(1-\frac{1}{q}\right)^{q+1}\right)q^2 = \left(1-\frac{1}{e}\right)q^2 +
O(q), \text{ as $q \to \infty$}.
\end{equation*}

\item[(b)] For $q$ sufficiently large, a Besicovitch set $B$
formed by the union of $q+1$ randomly chosen lines with distinct slopes will satisfy
\begin{equation*}
\left| \#B - \left(1-\frac{1}{e}\right)q^2 \right| < 2q \log q,
\end{equation*} with probability $1- O((\log q)^{-2})$.
In particular, the probability of randomly constructing a small Besicovitch
set tends to zero as $q \to \infty$.
\end{enumerate}
\end{thm}

As $1-1/e \approx 0.632$, we see that the average randomly chosen Besicovitch set will
contain around $0.632q^2$ points, whereas we expect a small Besicovitch set to consist
of about $0.5q^2$ points.

\section{Proofs of the results} \label{Proof Section}

\begin{proof}[Proof of the Incidence Formula]
Let us arbitrarily assign an ordering to the lines that comprise $B$: $\ell_0, \ldots,
\ell_q$. We use the fact that each pair of lines with distinct slopes must
intersect in exactly one point, and we argue essentially by inclusion--exclusion.

Fix $0 \leq j \leq q$. For $P$ a point on $\ell_j$, define $m_P(j)$ to be the number of
lines $\ell_i$ that contain $P$ with $i \leq j$. We wish to consider the intersections of
$\ell_j$ with $\ell_i$ for $i < j$. If all of these intersections are distinct, then
there are $q-j$ points on the line $\ell_j$ that do not lie on any $\ell_i$ with $i < j$.
For $P \in \ell_j$, we see $m_P(j) - 1$ of these lines meet at $P$; if $m_P(j) - 1 \geq 2$, then
we have undercounted the points on $\ell_j$ that do not lie on any $\ell_i$ with $i<j$
by $m_P(j)-2$ points. That is,
\[\#\left( \ell_j \setminus \bigcup_{i=0}^{j-1} \ell_i\right) = q- j + \sum_{P \in
\ell_j} \max\left\{0, m_P(j) - 2\right\}.\] Summing over all $j$ we get
\begin{equation*}
\begin{aligned}
\#B  &= \sum_{j=0}^q \#\left( \ell_j \setminus \bigcup_{i=0}^{j-1} \ell_i\right)
 = \sum_{j=0}^q (q-j) + \sum_{j=0}^q \sum_{P \in
\ell_j} \max\left\{0, m_P(j) - 2\right\} \\
&=  \frac{q(q+1)}{2} + \sum_{P \in B} \sum_{\substack{j=0 \\ P \in \ell_j}}^q
\max\left\{0,m_P(j) - 2\right\} \\
&=  \frac{q(q+1)}{2} + \sum_{P \in B} \max \left\{0, 1 + 2 + \cdots + \left(m_P -
2\right) \right\}  \\
&= \frac{q(q+1)}{2} + \sum_{P \in B} \frac{(m_P-1)(m_P-2)}{2}.
\end{aligned}
\end{equation*}
\end{proof}

\begin{example}Recall that we defined
\begin{equation*}
\label{Recall Minimal Set} B_0 = \left( \bigcup_{i \in \FF} \ell(i,-i^2) \right) \cup
\ell(\infty,0).
\end{equation*}

For $i,j$ distinct elements of $\FF$, one can easily see that
$\ell(i,-i^2)\cap \ell(j,-j^2) = \{(i+j,ij)\}$. Thus the lines
$\ell(i,-i^2), \ell(j,-j^2), \ell(k,-k^2)$ cannot share a common point
if $i,j,k$ are distinct. It follows that no point $P$ has multiplicity
$m_P > 3$, and if $m_P=3$, then $P$ must lie on the line
$\ell(\infty,0)$. In fact, $\ell(i,-i^2) \cap \ell(\infty,0) =
\{(0,-i^2)\}$.

If $i \not= 0$ and $q$ is odd, then precisely two of
our lines with nonzero slope pass through $(0,-i^2)$, namely $\ell(i, -i^2)$ and
$\ell(-i,-i^2)$. There are $\frac{q-1}{2}$ nonzero squares in $\FF$, so
\[
\sum_{P \in B_0} \frac{(m_P-1)(m_P-2)}{2} = \frac{q-1}{2}.
\]

If $q$ is even, then $\ell(i, -i^2)=\ell(-i,-i^2)$. There are no points of
multiplicity $m_P = 3$ in this case, and
\[
\sum_{P \in B_0} \frac{(m_P-1)(m_P-2)}{2} = 0.
\]

\end{example}

One can also prove the Incidence Formula in a fancier way using
intersection theory on algebraic surfaces. Roughly speaking, we
consider the divisor on $\mathbb{P}^2$ arising from $B$ consisting of
$q+1$ lines and compute its arithmetic genus in two ways: 1) using the
adjunction formula for divisors, and 2) by blowing up $\mathbb{P}^2$
at all of the multiple points of $B$ to get a surface on which the
lines in $B$ become pairwise disjoint. We leave the details to the
interested reader. (See \cite[Exercise V.1.3 and Corollary V.3.7]{Hart}.)

In order to prove Theorems~\ref{Unconditional Theorem} and
\ref{Conditional Theorem}, we require the following lemma:

\begin{tpl} \label{Triple Point Lemma}
Let $q$ be odd. Suppose $B$ is a Besicovitch set with $B=\bigcup_{i \in \FF \cup \{\infty\}}
\ell(i,b_i)$. Then with at most one exception, for any choice of $i \in \FF\cup \{\infty\}$ there
exists a point $P \in \ell(i, b_i)$ with $m_P \geq 3$.
\end{tpl}

\begin{proof}
  Suppose $\ell$ and $\ell'$ are two lines in $B$ such that no point
  $P$ with $m_P \geq 3$ lies on either one. Without loss of
  generality, we may apply a translation followed by a linear
  automorphism of $\FF^2$ so that it suffices to assume $\ell(0,0)$
  and $\ell(\infty,0)$ are the two lines. Note that translations and
  linearautomorphisms carry lines to lines and respect the multiplicities $m_P$.

As $i$ varies through $\FF^{\times}$, it must be true that the $y$-intercepts of
$\ell(i,b_i)$ are distinct. For if not, there would exist a triple point on the line
$\ell(\infty,0)$. Similarly, the $x$-intercepts of these lines must be distinct. Note
that $b_i$ cannot be zero for any $i \not=0$ since that would imply the existence of a
triple point at the origin. The $x$- and $y$-intercepts of $\ell(i,b_i)$ are $-i/b_i$ and
$b_i$, respectively. We conclude that \[\{i: i \in \FF^{\times}\} = \{-i/b_i: i \in
\FF^{\times}\} = \{b_i: i \in \FF^{\times}\},\] since each set is a collection of $q-1$
distinct nonzero elements of $\FF$. Using the fact that
the product of all nonzero elements of $\FF$ is $-1$ when $q$ is odd, we see that
\begin{equation*}
-1 = \prod_{i \in \FF^{\times}} i = \prod_{i \in \FF^{\times}}
\left(-\frac{i}{b_i}\right) = (-1)^{q-1} \frac{\prod_{i \in \FF^{\times}} i}{\prod_{i \in
\FF^{\times}} b_i} = 1.
\end{equation*}
Evidently this is a contradiction, so we are forced to accept the statement of the lemma.
\end{proof}

\begin{proof}[Proof of Theorem~\ref{Unconditional Theorem}]
We may suppose that $B$ consists of $q+1$ lines, arbitrarily labelled $\ell_0, \ldots,
\ell_q$. For a point $P \in B$, there are $m_P$ lines passing through it; we make the
trivial observation
\begin{equation*}
\frac{(m_P-1)(m_P-2)}{2} = \frac{1}{m_P}\sum_{\substack{j=0 \\ P \in \ell_j}}^q
\frac{(m_P-1)(m_P-2)}{2}.
\end{equation*}
It follows that
\begin{equation} \label{Triple Point Inequality}
\begin{aligned}
\sum_{P \in B} \frac{(m_P-1)(m_P-2)}{2} &= \sum_{P \in B} \sum_{\substack{j=0 \\ P \in
\ell_j}}^q \frac{m_P^2 - 3m_P + 2}{2m_P} \\
&= \sum_{j=0}^q \sum_{P \in \ell_j} \frac{m_P^2 - 3m_P + 2}{2m_P} \\
&\geq \sum_{j=0}^q
\sum_{\substack{P \in \ell_j \\ m_P \geq 3}} \frac{1}{3}.
\end{aligned}
\end{equation}
For the inequality, note that the function $x \mapsto \frac{x^2-3x+2}{2x}$ is increasing for $x \geq 3$
and evaluates to $1/3$ for $x=3$. By the Triple Point Lemma, we know that every line,
except perhaps one, contains a point of multiplicity three or greater. Hence there are at
least $q$ terms in the final double sum in \eqref{Triple Point Inequality}, and we obtain
\begin{equation*}
\sum_{P \in B} \frac{(m_P-1)(m_P-2)}{2} \geq \frac{q}{3}.
\end{equation*}

\end{proof}

\begin{proof}[Proof of Theorem~\ref{Conditional Theorem}]
Without loss of generality, we may apply a linear automorphism of $\FF^2$ and assume that
all points $P \in B$ with $m_P \geq 3$ lie on the line $\ell(\infty,0)$. Suppose the
number of such points is $T$. Let us agree to write $\sum'$ for the sum over
points $P$ with $m_P \geq 3$. As every line $\ell(i,b_i)$ with $i \not= \infty$ must intersect
$\ell(\infty,0)$ exactly once, and there exists at most one line that does not contain a point $P$
with $m_P \geq 3$, we find that the sum of the multiplicities $m_P$ over all points with
$m_P \geq 3$ must satisfy $\sum' m_P = q + T - \delta$, where $\delta$ is $0$~or~$1$.  We
now have
\begin{equation} \label{First estimate}
\begin{aligned}
\sum_{P \in B} \frac{(m_P-1)(m_P-2)}{2}
&= \frac{1}{2}\primesum m_P^2 - \frac{3}{2}\primesum m_p + \primesum 1 \\
&= \frac{1}{2}\primesum m_P^2 - \frac{3}{2}(q -\delta) -\frac{1}{2}T.
\end{aligned}
\end{equation}
By the Cauchy-Schwartz inequality, we find that
\begin{equation} \label{Second estimate}
\left(\primesum m_p\right)^2 \leq \left(\primesum 1 \right) \left(\primesum m_P^2 \right)
= T \primesum m_P^2.
\end{equation}
Combining \eqref{First estimate} and \eqref{Second estimate}, and again using the fact
that $\sum' m_P = q+T-\delta$, we find that
\begin{equation*}
\begin{aligned}
\sum_{P \in B} \frac{(m_P-1)(m_P-2)}{2} &\geq \frac{1}{2T} (q+T-\delta)^2
-\frac{3}{2}(q-\delta) -\frac{1}{2}T \\ &= \frac{1}{2T}\left(q-\delta\right)^2 -
\frac{1}{2}(q-\delta).
\end{aligned}
\end{equation*}
This last expression is a decreasing function of $T$. At least two non-vertical lines
(slope $i \not= \infty$) pass through each point $P$ with $m_P \geq 3$, and at most $q$
non-vertical lines pass through these points in total. So $2T \leq q$, but since $q/2$ is
not an integer, we obtain  $T \leq (q-1)/2$.  We now see that
\begin{equation} \label{Final expression}
\begin{aligned}
\sum_{P \in B} \frac{(m_P-1)(m_P-2)}{2} &\geq
\frac{\left(q-\delta\right)^2}{q-1} - \frac{1}{2}(q-\delta) \\
&= \frac{q-1}{2} + \frac{3}{2} - \frac{3}{2}\delta + \frac{(1-\delta)^2}{q-1}.
\end{aligned}
\end{equation}
The final three terms contribute a non-negative quantity for $\delta=0$ or $\delta=1$,
which shows the desired inequality.

As for the final claim of the theorem, equality clearly holds if
$m_P=3$ for $(q-1)/2$ points $P \in B$ and $m_P <3$ otherwise.
Conversely, if equality holds in the theorem, then equality must hold
in \eqref{Final expression}.  Evidently this is equivalent to saying
$T=(q-1)/2$. Now there are $(q-1)/2$ nonzero terms in the sum in
\eqref{Final expression}, and their sum must be $(q-1)/2$. We conclude
that $\frac{(m_P-1)(m_P-2)}{2}=1$ for all $P$ with $m_P \geq 3$.  That
is, $m_P=3$ for exactly $(q-1)/2$ points, and $m_P<3$ for all other
points in $B$.
\end{proof}

To prove Theorem~\ref{Probability Theorem}, we first formalize the underlying probability space.
Let $\Omega = \bigoplus_{i \in \FF \cup \{\infty\}} \FF$. We can
identify an element $\sum b_i \in \Omega$ with a Besicovitch set by setting $B=
\bigcup_{i \in \FF \cup \{\infty\}} \ell(i, b_i)$. We will use this identification without further
comment. We make $\Omega$ into a probability
space by assigning probability $q^{-(q+1)}$ to each Besicovitch set.

\begin{proof}[Proof of Theorem~\ref{Probability Theorem}]
We will proceed in three steps. The first is to calculate the mathematical expectation
for the cardinality function $\#: \Omega \to \RR$.

For $P \in \FF^2$, let $f_P:\Omega \to \RR$ be the characteristic function of $P$; i.e.,
\[f_P(B) =
\begin{cases}
1,& \text{if } P \in B, \\
0,& \text{otherwise.}
\end{cases}
\]
It follows that $\#B = \sum_{P \in \FF^2} f_P(B)$.

For a given point $P \in \FF^2$ we now calculate $\prob\{f_P=1\}$, the probability that $P$
appears in a randomly chosen Besicovitch set. For fixed $i \in \FF \cup
\{\infty\}$, the probability that $P$ does not lie on $\ell(i, b_i)$ is $1-1/q$,
since there are $q$ choices for the $y$-intercept $b_i$. The probability that $P$ does
not lie in a given Besicovitch set $B \in \Omega$ is the probability that it lies on none
of the lines comprising B. As the $y$-intercepts of lines with distinct slopes are
independent, we see that
\begin{equation} \label{Probability of zero}
\prob\{f_P=0\} = \prod_{i \in \FF
\cup \{\infty\}} \left(1-\frac{1}{q}\right) = \left(1-\frac{1}{q}\right)^{q+1}.
\end{equation}
Hence
\begin{equation} \label{Single Probability}
\prob\{f_P=1\} = 1-\left(1-\frac{1}{q}\right)^{q+1}.
\end{equation}

We can now determine the expectation of the cardinality function:
\begin{equation} \label{Expectation}
\begin{aligned}
\EE(\#) &= \sum_{B \in \Omega} \#B \cdot \prob \{B\}= \sum_{B \in \Omega} \sum_{P \in
\FF^2} f_P(B) \cdot \prob \{B\} \\
&= \sum_{P \in \FF^2} \sum_{B \in \Omega} f_P(B) \cdot \prob \{B\}  = \sum_{P \in \FF^2}
\prob\{f_P=1\} = \left(1- \left(1-\frac{1}{q}\right)^{q+1}\right)q^2 \\ &=
\left(1-\frac{1}{e}\right)q^2 + O(q), \text{ as $q \to \infty$}.
\end{aligned}
\end{equation}
This completes part (a) of the theorem.

The second step in the proof is to compute the variance of the cardinality function. To
this end, we will need to determine $\prob\{f_P=f_Q=1\}$ for two distinct points $P, Q
\in \FF^2$. We can rewrite this probability as
\begin{equation} \label{rephrase probability}
\begin{aligned}
\prob\{f_P=f_Q=1\} &= 1 - \prob\{f_P=f_Q=0\} \\
& \qquad \quad - \prob\{f_P=1,f_Q=0\} - \prob\{f_P=0,f_Q=1\} \\
&= 1 + \prob\{f_P=f_Q=0\} - \prob\{f_P=0\}- \prob\{f_Q=0\}.
\end{aligned}
\end{equation}
The second term is the only one we don't know yet. There is precisely one line containing
both $P$ and $Q$, say $\ell(j,a)$. The probability that a line with slope $j$ doesn't
contain $P$ or $Q$ must be $1-1/q$. For any other slope $i \not= j$, there is
precisely one line with slope $i$ passing through $P$, and one through $Q$. The
probability that a line with slope $i \not= j$ does not contain $P$ or $Q$ is
$1-2/q$. Again by independence of $y$-intercepts it follows that
$\prob\{f_P=f_{Q}=0\} = (1-1/q)(1-2/q)^{q}$. We
conclude from \eqref{rephrase probability} and \eqref{Probability of zero} that
\begin{equation} \label{joint probability}
\prob\{f_P=f_Q=1\} = 1 + \left(1-\frac{1}{q}\right)\left(1-\frac{2}{q}\right)^{q} -
2\left(1-\frac{1}{q}\right)^{q+1}.
\end{equation}

The variance of the cardinality function is given by
\begin{equation*} \label{variance}
\begin{aligned}
\Var(\#) &= \EE(\#^2) - \EE(\#)^2 = \sum_{B \in \Omega} \sum_{P,Q \in \FF^2} f_P(B)f_Q(B)
\cdot \prob\{B\} - \EE(\#)^2\\
&= \sum_{P,Q \in \FF^2} \prob\{f_P=f_Q=1\} - \EE(\#)^2\\
&= \sum_{P \not= Q \in \FF^2}\prob\{f_P=f_Q=1\} + \sum_{P\in
\FF^2} \prob\{f_P=1\} - \EE(\#)^2\\
&= q(q+1)(q-1)^2\left(1-\frac{2}{q}\right)^q +
q^2\left(1-\frac{1}{q}\right)^{q+1}\left\{1-q^2\left(1-\frac{1}{q}\right)^{q+1}\right\}
\\ &= \left(\frac{1}{e} -\frac{5}{2e^2}\right)q^2 + O(q), \text{ as $q \to \infty$}.
\end{aligned}
\end{equation*}
The second to last step follows from \eqref{joint probability}, \eqref{Single
Probability}, \eqref{Expectation}, and a bit of simplification.

The third and final step of the proof is an application of the Chebyshev inequality.
Recall that the Chebyshev inequality asserts that for any function $g:\Omega \to \RR$ and
any $\varepsilon > 0$, we have
\[\prob\{B \in \Omega:|g(B) - \EE(g)| \geq \varepsilon\} \leq \frac{1}{\varepsilon^2} \Var(g).\]
Applying this to our situation with $g = \#$ and $\varepsilon = q\log q$ shows
\[\prob\{B \in \Omega: \left| \#B - \EE(\#) \right| \geq q \log q \}  = O((\log
q)^{-2}).\] As $\EE(\#)$ differs from $\left(1-1/e\right)q^2$ by $O(q)$,
part (b) of the theorem follows.
\end{proof}

\section{Acknowledgments}

I would like to thank Andrew Granville for presenting this problem to me, and for all of the
helpful conversations we've had on this subject. Eric Pine and Boris Mezhericher provided
much of the computational evidence for Conjectures~\ref{Size Conjecture}~and~\ref{Multiplicity Conjecture}.
This work was primarily supported by the NSF VIGRE grant at the University of Georgia.

\end{document}